\documentclass[11pt]{amsart}

\usepackage{amssymb, amsthm, amsmath, gensymb}
\usepackage{graphicx, comment, accents}
\usepackage[small]{caption}
\usepackage{subcaption}
\usepackage{epsfig}
\usepackage{tikz, pgfplots, float}
\usepackage{stmaryrd}

\newcommand{\later}[1]{}
\newcommand{\old}[1]{}

\usepackage{amsfonts}

\usepackage[utf8]{inputenc}
\usepackage{fullpage}
\usepackage{framed}
\usepackage{multirow}
\usepackage{enumerate}
\usepackage{url}
\usepackage[breaklinks]{hyperref}
\usepackage{cleveref}
\hypersetup{
	colorlinks = true, 
	urlcolor = cyan, 
	linkcolor = teal, 
	citecolor = cyan 
}

\newcommand{\ex}{\mathop{}\!\mathrm{ex^*}}
\newcommand{\orig}{\mathop{}\!\mathrm{orig}}
\newcommand{\dist}{\mathop{}\!\mathrm{dist}}
\newcommand{\re}{\mathop{}\!\mathrm{re^*}}
\newcommand{\LA}{\mathop{}\!\mathrm{La^*}}
\newcommand{\LAw}{\mathop{}\!\mathrm{La}}
\newcommand{\RLA}{\mathop{}\!\mathrm{La_R^*}}
\newcommand{\RLAw}{\mathop{}\!\mathrm{La_R}}
\newtheorem{thm}{Theorem}[section]

\newtheorem{lem}[thm]{Lemma}
\newtheorem{cor}[thm]{Corollary}
\newtheorem{claim}[thm]{Claim}
\newtheorem{cons}[thm]{Construction}
\newtheorem{prop}[thm]{Proposition}
\theoremstyle{definition}
\newtheorem{rem}[thm]{Remark}
\newtheorem{prob}[thm]{Problem}

\newcommand{\cC}{{\mathcal C}}
\newcommand{\cP}{{\mathcal P}}
\newcommand{\cF}{{\mathcal F}}

\newcommand{\cG}{{\mathcal G}}

\newcommand{\cA}{{\mathcal A}}

\newcommand{\cT}{{\mathcal T}}

\title{Rainbow Tur\'an problems for forbidden subposets}
\author{Bal\'azs Patk\'os}
\address{HUN-REN Alfr\'ed R\'enyi Institute of Mathematics} 
\email{patkos@renyi.hu}
\date{}

\begin{document}

\begin{abstract}
    A family $\cG$ of sets is a copy of a poset $(P,\leqslant)$ if $(\cG,\subseteq)$ is isomorphic to $(P,\leqslant)$. The forbidden subposet problem asks for determining $\LA(n,P)$, the maximum size of a family $\cF\subseteq 2^{[n]}$ that does not contain any copy of $P$. We study the rainbow version of this problem: what is the maximum size $\RLA(n,P)$ of a family $\cF=\cup_{i=1}^mA^i$ such that all $A^i$ are antichains and there is no copy of $P$ with all sets coming from distinct $A^i$ or equivalently $\cF$ admits a proper coloring (sets $F\subset F'$ must receive different colors) with no rainbow copy of $P$.

    A poset $(Q,\leqslant')$ rainbow forces $(P,\leqslant)$ if any proper coloring $c$ of $Q$ ($q\leqslant' q'$ or $q'\leqslant' q$ implies $c(q)\neq c(q')$) admits a rainbow copy of $P$. We establish connection between the $\LA$ and the $\RLA$ functions via poset rainbow forcing, determine the asymptotics of  $\RLA(n,T)$ for all tree posets and obtain further exact or asymptotic results for antichains and complete bipartite posets.
\end{abstract}

\maketitle

\section{Introduction}

In extremal combinatorics, most of the problems ask for determining the largest, smallest or otherwise extremal combinatorial structure that possesses some prescribed property. In this paper, we will investigate such problems about families of finite sets. We use standard notation: $[n]$ stands for the set of the first $n$ positive integers, for any set $X$, the notations $2^X$, $\binom{X}{k}$ are used for the power set of $X$ and the family of all $k$-element subsets of $X$, respectively. A subfamily $\cG\subseteq \cF$ is a copy\footnote{In the literature of forbidden subposet problems, this kind of copies are called \textit{strong} or \textit{induced} copies of $P$. Since we will only consider other types of copies in the concluding remarks section, we will omit the adjective.} of a poset $(P,\leqslant)$ if $\cG$ is isomorphic to $(P, \leqslant)$ when ordered by inclusion. By abuse of notation, we will write $P$ instead of $(P,\leqslant)$. A family $\cF$ is $P$-free if it does not contain any copies of $P$. For a set $\cP$ of posets, $\cF$ is $\cP$-free if it is $P$-free for all $P\in \cP$. The problem of determining the maximum size $\LA(n,P)$ ($\LA(n,\cP)$) of a $P$-free ($\cP$-free) family $\cF\subseteq 2^{[n]}$ (usually referred to as the \textit{forbidden subsposet problem}) was introduced by Katona and Tarj\'an \cite{KT} as a natural generalization of theorems of Sperner \cite{S} and Erd\H os \cite{E} considering the cases $P=C_2$ and $P=C_k$, respectively, where $C_k$ denotes the chain (totally ordered set) on $k$ elements. For relatively recent surveys, we refer to \cite{GP,GLsurv}. 

In a coloring of a combinatorial object, a substructure is \textit{rainbow} if all its elements receive distinct colors. In the case of rainbow extremal problems, a coloring is assumed to exist, and the size of the largest substructure without rainbow forbidden object is addressed. If all elements are assigned the same color, then we cannot expect any rainbow object. We say that a coloring of a poset $P$ is \textit{proper} if $c(p)=c(q)$ implies $p\perp q$, i.e. $p$ and $q$ are incomparable. The \textit{comparability graph} $CG(P)$ of a poset $P$ has vertex set $P$ with $xy$ being an edge if and only if $x$ and $y$ are comparable. So a proper coloring of $P$ is equivalent to a proper coloring of its comparability graph. The general rainbow Tur\'an problem is as follows: $P'\subset Q$ for $(Q,\leqslant')$ is a copy of a poset $(P,\leqslant)$ if $(P,\leqslant)$ and $(P',\leqslant')$ are isomorphic, $Q$ is $P$-free if it does not contain any copies of $P$. For posets $P$ and $Q$, determine the maximum size $\RLA(Q,P)$ of $A\subset Q$ such that $Q[A]$, the subposet of $Q$ induced by $A$ can be colored properly without admitting a rainbow copy of $P$. Equivalently, what is the maximum size of a union of antichains in $Q$ such that there is no copy of $P$ with all elements coming from different antichains. We will be interested in the case $Q=B_n=(2^{[n]},\subseteq)$ and write $\RLA(n,P)$ for $\RLA(B_n,P)$.

There is a very natural way how to obtain large families of sets without forbidden copies of $P$. Let $e^*(P)$ denote the maximum positive integer $k$ such that the union of $k$ middle layers of $B_n$ is $P$-free for all $n$. By definition, $\LA(n,P)\ge (e^*(P)+o(1))\binom{n}{\lfloor n/2\rfloor}$ and a longstanding conjecture \cite{B,GLu,P}  stated equality should hold for all posets $P$. This was refuted by Ellis, Ivan, and Leader \cite{EIL}. Similarly to $e^*(P)$, let $\re(P)$ denote the maximum positive integer $k$ such that the $k$ middle layers of $B_n$ do not contain a copy of $P$ with all elements taken from distinct layers. By definition, we have the following statement.

\begin{prop}
    $\RLA(n,P)\ge (\re(P)+o(1))\binom{n}{\lfloor n/2\rfloor}$.
\end{prop}

The example of \cite{EIL} showing the existence of a poset $P$ with $\liminf_{n\rightarrow \infty} \frac{\LA(n,P)}{\binom{n}{\lfloor n/2\rfloor}}>e^*(P)$ does not yield an example for a poset with $\liminf_{n\rightarrow \infty} \frac{\RLA(n,P)}{\binom{n}{\lfloor n/2\rfloor}}>\re(P)$. In this paper, we will determine $\RLA(n,P)$ for several classes of posets, but none of them will possess this property.

\begin{prob}
    Find a poset $P$ with $\liminf_{n\rightarrow \infty} \frac{\RLA(n,P)}{\binom{n}{\lfloor n/2\rfloor}}>\re(P)$ (or prove that there does not exist any).
\end{prob}

In our investigations, the following notion will be crucial. We say that a poset $Q$ \textit{rainbow forces} a poset $P$ if in any proper coloring of $Q$ there exists a rainbow copy of $P$. For any poset $P$, we set $F(P)$ to be the set of posets that rainbow force $P$ and write $M(P)$ to denote the inclusion-wise minimal elements of $F(P)$, i.e. those posets $Q$ that rainbow force $P$, but for any $q\in Q$, the poset $Q\setminus \{q\}$ does not rainbow force $P$. The connection between ordinary forbidden subposet problems and rainbow Tur\'an problems is established in the following observation.

\begin{prop}\label{connection}
    For any poset $P$, we have $\RLA(n,P)=\LA(n,M(P))$.
\end{prop}

A theorem of Methuku and P\'alv\"olgyi \cite{MP} states that $\LA(n,P)=O(\binom{n}{\lfloor n/2\rfloor})$ for any poset $P$. This, by Proposition \ref{connection}, implies $\RLA(n,P)=O(\binom{n}{\lfloor n/2\rfloor})$ for any poset $P$. The family $\binom{[n]}{\lfloor n/2\rfloor}$ shows that this order of magnitude is achievable without having a rainbow copy of $P$ (all sets are colored the same) if $P$ contains at least 2 elements. All our upper bounds on $\RLA(n,P)$ will follow by upper bounds on $\LA(n,Q)$ for some $Q\in F(P)$. Heuristically, one may think that 'smaller' posets are harder to avoid than 'larger' ones.  
The minimum size of a poset in $M(P)$ is denoted by $m(P)$. We have the following general statement.
\begin{thm}\label{mP}
    For any poset of size $k$, we have $k \le m(P) \le \binom{k+1}{2}$. Furthermore, $m(P)=k$ if and only if $P=C_k$, and $m(P)=\binom{k+1}{2}$ if and only if $P=A_k$, the antichain on $k$ elements.
\end{thm}

The \textit{directed Hasse diagram} $\overrightarrow{H}(P)$ of a poset $P$ is the directed graph with vertex set $P$ and $(xy)\in E(H(P))$ if $x\leqslant_P y$ and there is no $z\in P$ with $z\neq x,y$ $x\leqslant_P z\leqslant_P y$, in other words $y$ \textit{covers} $x$. If one removes the orientations of the arcs of $\overrightarrow{H}(P)$, then one obtains the (undirected) Hasse diagram of $P$. A poset is a \textit{tree poset} if its undirected Hasse diagram is a tree. For any poset $P$, its \textit{height} $h(P)$ is the length of its longest chain. An important class of posets for which the asymptotics of $\LA(n,P)$ is known is that of tree posets, and this asymptotic value depends only on the height of the tree poset $T$.

\begin{thm}[Boehnlein, Jiang \cite{BJ}]\label{bj}
    For any tree poset $T$, we have $\LA(n,T)=(h(T)-1+o(1))\binom{n}{\lfloor n/2\rfloor}$.
\end{thm}

One of the main results of this paper determines the asymptotics of $\RLA(n,T)$ if $T$ is a tree poset. By showing the existence of a universal tree poset $\cT^k$ of height $k$ that rainbow forces every tree poset on $k$ elements and applying Proposition \ref{connection} and Theorem \ref{bj}, we shall obtain the following result.

\begin{thm}\label{tree}
    For any tree poset $T$, we have $\RLA(n,T)=(|T|-1+o(1))\binom{n}{\lfloor n/2\rfloor}$.
\end{thm}

The complete $r$-level poset $K_{s_1,s_2,\dots,s_r}$ has elements $x^1_1,\dots,x^1_{s_1},x^2_1,\dots,x^2_{s_2},\dots,x^r_1,\dots,x^r_{s_r}$ with the only relations being $x^j_i<x^{j'}_{i'}$ whenever $j<j'$ for all $1\le i\le s_j$, $1\le i'\le s_{j'}$. As $K_{s,t}$ is a subposet of $K_{s,1,t}$, and $K_{s,1,t}$ is a tree poset, Theorem \ref{bj} implies $\LA(n,K_{s,t})\le(2+o(1))\binom{n}{\lfloor n/2\rfloor}$. The family of the middle two levels of $B_n$ shows that the inequality is an equality. $\Sigma(n,k)$ denotes the sum $\sum_{i=1}^k\binom{n}{\lfloor\frac{n-k}{2}\rfloor+i}$ of the $k$ largest binomial coefficients of order $n$, the size of the family consisting of the $k$ largest layers of $2^{[n]}$.

\begin{prop}\label{easylar}
\begin{enumerate}[(i)]
    \item 
    For any $k\ge 2$, there exists $n_0=n_0(k)$ such that if $n\ge n_0$, then $\RLA(n,A_k)=\Sigma(n,k-1)+2$.
    \item 
    For all $s,t\ge 2$, we have $\RLA(n,K_{s,t})=(s+t+o(1))\binom{n}{\lfloor n/2\rfloor}$.
\end{enumerate}
\end{prop}

The smallest poset for which the asymptotics of $\LA(n,P)$ is yet to be determined is the diamond poset $\Diamond=K_{1,2,1}$. We obtain the following bounds on $\RLA(n,\Diamond)$.

\begin{thm}\label{diamond}
    $(3+o(1))\binom{n}{\lfloor n/2\rfloor}\le \RLA(n,\Diamond)\le (3+\frac{1}{6}+o(1))\binom{n}{\lfloor n/2\rfloor}$.
\end{thm}

\begin{rem}
    Both of the main concepts, $\RLA(n,P)$ and rainbow forcing can be defined for graphs and directed graphs, too. One can write $\ex(G,F)$ to denote the size of the largest vertex subset $U$ of $V(G)$ such that $G[U]$ can be properly colored without a rainbow copy of $F$. A graph $H$ rainbow forces $F$ if any proper coloring of $H$ contains a rainbow copy of $F$. Some of our proofs will work also in this context. We will comment on the analogies and differences after most proofs in the paper.
\end{rem}

The rest of the paper is organized as follows: Section 2 contains basic results on rainbow forcing posets, in particular, the proof of Theorem \ref{mP}. Then in Section 3, we prove results on $\RLA(n,P)$, in particular Theorem \ref{tree} and Theorem \ref{diamond}. Section 4 contains some concluding remarks.

\medskip

\textbf{Notation and terminology.} A poset $P$ can be canonically partitioned into antichains by taking $A^1$ to be the set of $P$'s minimal elements, then $A^j$ to be the set of minimal elements of $P\setminus \cup_{i=1}^{j-1}A^i$. Then the \textit{rank} of $x\in P$ is $r(x)=j$ with $x\in A^j$. The \textit{dual} of a poset $P=(P,\leqslant)$ is $P'=(P,\leqslant')$ with $x\leqslant' y$ if and only if $y\leqslant x$. Finally, for any poset $P$ and elements $x,z\in P$, we write $U_P(x)=\{y\in P:x\leqslant y\}$ and $D_P(x)=\{y\in P:y\leqslant x\}$, and $[x,z]=U_P(x)\cap D_P(z)$.

\section{Results on poset rainbow forcing}

In this section, we prove Theorem \ref{mP} and some basic results about poset rainbow forcing. Theorem \ref{mP} will immediately follow from the next lemma. We need a definition first. Let $P$ be a poset with elements $p_1,p_2,\dots,p_k$ and let $\pi$ be a permutation of $[k]$. Then  let $\perp_P(\pi)=k+\sum_{j=1}^k|\{i<j:p_{\pi(i)}\perp p_{\pi(j)}\}|$. Observe that for any $\pi,\pi'$, we have $\perp_P(\pi)=\perp_P(\pi')=k+\binom{k}{2}-|E(CG(P))|$ and we write $\perp(P)$ for this value. 

\begin{lem}\label{perp}
    For any poset $P$, we have $m(P)\le \perp(P)$.
\end{lem}

\begin{proof}
    Let $P$ be a poset on $k$ elements $p_1,p_2,\dots,p_k$ and let us fix a permutation $\pi$ of $[k]$. Consider the following blow-up $P^\pi$ of $P$: $p_{\pi(j)}$ is replaced with a chain of length $1+|\{i<j:p_{\pi(i)}\perp p_{\pi(j)}\}|$ and two elements in the chains of $p_{\pi(i)}$ and $p_{\pi(j)}$ inherit the relationship of $p_{\pi(i)}$ and $p_{\pi(j)}$ in $P$. We claim that $P^\pi$ rainbow forces $P$ for any permutation $\pi$. First observe that if one takes one element from each chain, then they form a copy of $P$. The elements of a chain receive different colors in any proper coloring of $P^\pi$. Also, the elements of two chains corresponding to comparable elements of $P$ must receive different colors in any proper coloring $c$ of $P^\pi$. So for any proper coloring of $P^{\pi}$, we can pick elements from each chain greedily in the order of $\pi$: as the chain of $p_{\pi(j)}$ consists of $1+|\{i<j:p_{\pi(i)}\perp p_{\pi(j)}\}|$ elements, one of them has to receive a color different from all colors of previously picked vertices from incomparable chains, and their color must be different from previously picked elements from comparable chains as the coloring is proper. So there is a rainbow copy of $P$ with respect to $c$.
\end{proof}

\begin{proof}[Proof of Theorem \ref{mP}]
    If $Q$ forces $P$, then $P$ must be a subposet of $Q$ so $|P|\le |Q|$. If $|P|=|Q|$, then every proper coloring of $Q$ must be rainbow, but that is true only for chains. Also, for this reason $C_k$ rainbow forces itself for any $k$.

    As $1+|\{i<j:p_{\pi(i)}\perp p_{\pi(j)}|\le j$ for any poset and permutation, we obtain $m(P)\le \perp(P)\le \sum_{j=1}^kj=\binom{k+1}{2}$. Also, if $P\neq A_k$, then for some $i<j$ we have $p_i\not\perp p_j$ and so we have strict inequality above. Finally, we need to show $m(A_k)=\binom{k+1}{2}$. We proceed by induction on $k$ with the case $k=1$ being trivial as $A_1=C_1$. Consider a poset $P\in M(A_k)$. We claim that it has to contain a chain of length $k$. Indeed, if not then by Mirsky's theorem \cite{M}, $P$ is the union of at most $k-1$ antichains, and thus it can be properly colored using at most $k-1$ colors, and so we would not have a rainbow copy of $A_k$ (or any rainbow poset on $k$ elements). Next we claim that if we remove the elements of this $k$-chain $C$ from $P$, then the remaining poset $P^-$ rainbow forces $A_{k-1}$. Indeed, if not, then there exists a proper coloring $c$ of $P^-$ without a rainbow copy of $A_{k-1}$. But then extending $c$ to a proper coloring of $P$ (this is possible by using completely new colors), we obtain a proper coloring that has no rainbow $A_k$ in $P$ as we can pick at most one element of $C$ into an antichain. So by induction, we obtain $|P|=k+|P^-|\ge k+\binom{k}{2}=\binom{k+1}{2}$ as claimed.
\end{proof}

\begin{rem}
    The proof of Lemma \ref{perp} and almost all proof of Theorem \ref{mP} work in the case of graphs if we replace incomparability with non-adjacency. So the minimum order of a graph $G$ that rainbow forces $F$ is between $|V(F)|$ and $\binom{|V(F)|+1}{2}$ and $m(F)=|V(F)|$ if and only if $F$ is a complete graph. Moreover, if $F$ is not an empty graph, then $m(F)<\binom{|V(F)|+1}{2}$.
\end{rem}

For posets $P$ and $Q$, we write $P\olessthan Q$ to denote the poset on $P\cup Q$ with relations inherited from $P$ and $Q$ and with relations $p\leqslant q$ for any $p\in P,q\in Q$.
\begin{lem}\label{linsum}
    For any poset $P_1,P_2$, if $Q_1\in M(P_1)$ and $Q_2\in M(P_2)$, then $Q_1'\in M(P_1')$ and $Q_1\olessthan Q_2\in M(P_1\olessthan P_2)$. 
\end{lem}

\begin{proof}
    The statement about the dual is immediate. As the sets of colors used for elements of $Q_1$ and $Q_2$ are disjoint in any proper coloring of $Q_1\olessthan  Q_2$, $Q_1\olessthan Q_2$ rainbow forces $P_1\olessthan P_2$. For any $x\in Q_2$, we need to show $Q_1\olessthan Q_2\setminus \{x\}$ does not rainbow force $P_1\olessthan P_2$ (the case $x\in Q_1$ is analogous). As $Q_2\in M(P_2)$ there exists a proper coloring $c_2$ of $Q_2\setminus \{x\}$ that does not admit a rainbow $P_2$. Let $z$ be an arbitrary maximal element of $Q_1$. As $Q_1\in M(P_1)$ there exists a proper coloring $c_1$ of $Q_1\setminus \{z\}$ that does not admit a rainbow $Q_1$. By recoloring but keeping the color classes of $c_1$, we can assume that the sets of colors used by $c_1$ and $c_2$ are disjoint. Let $c$ be the coloring of $Q_1^{Q_2}\setminus \{x\}$ defined by $c(y)=c_2(y)$ if $y\in Q_2\setminus \{x\}$, $c(y)=c_1(y)$ if $y\in Q_1\setminus \{z\}$ and by letting $c(z)$ be a completely new color. 
    
    Observe that $c$ is proper as colors used in $c_1$ are not used in $c_2$ and vice versa. We claim that $c$ does not admit a rainbow copy of $P_1\olessthan P_2$. Suppose it does. Then the $P_2$ part must contain an element $p$ from the $Q_1$ part as $c_2$ does not admit a rainbow $P_2$. Also, as $c_1$ does not admit a rainbow $P_1$, the $P_1$ part of the rainbow copy of $c$ must contain either $z$ or an element $p'$ from $Q_2$. As $p$ is in the $P_2$ part, $p'$ or $z$ is in the $P_1$ part, $p'$ or $z$ should be smaller than $p$, but  $p',z\not\leqslant_{Q_1\olessthan Q_2} p$ as $p' \in Q_1$, and either $z$ is a maximal element of $Q_1$ or $p'\in Q_2$.
\end{proof}

Let $O_k$ denote the \textit{organ} $\bigoplus_{i=1}^kC_i$, where $\bigoplus_{i=1}^kP_i$ denotes the pairwise incomparable disjoint union of the posets $P_i$. Let $H_k$ denote the \textit{harp} $(C_1\olessthan O_k)\olessthan C_1$, i.e. $O_k$ with a smallest and a largest element added. We introduce the posets $\vee=K_{1,2}$ and $J=C_1\olessthan O_2$.

\begin{lem}\label{smallexample}
    $M(A_2)=\{O_2\}$, $M(\vee)=\{J\}$, $H_2\in M(\Diamond)$, $O_k\in M(A_k)$.
\end{lem}

\begin{proof}
$O_2$ rainbow forces $A_2$ as in any proper coloring, at least one element of $C_2$ receives different color other than the element in $C_1$. $O_2\in M(A_2)$ as removing any element would either remove all $A_2$s or would leave $A_2$ that is 1-colorable. Lemma \ref{linsum} implies $J \in M(\vee), H_2\in M(\Diamond)$.

If a poset $P$ is $O_2$-free, then it is a complete multipartite poset. This can be seen by induction on $|P|$: antichains are special complete multipartite posets with one part. Otherwise consider the set $M$ of all minimal elements of $P$. As $P$ is $O_2$-free, $P=M\olessthan (P\setminus M)$, and one can apply induction to $P\setminus M$. Complete multipartite posets are properly colorable by coloring according to parts. This coloring does not contain a rainbow $A_2$, so $O_2$ is the only poset that minimally rainbow forces $A_2$.

If $P\in M(\vee)$ has a smallest element $x$, then $P\setminus \{x\}$ minimally rainbow forces $A_2$ and so $P=J$. If $P$ has multiple minimal elements $x_1,x_2,\dots,x_s$, then if $P$ does not contain $J$, then all $U_P(x_i)$s are complete multipartite posets. So for any $i$ the upset $U_P(x_i)$ is of the form $\cup_{j=1}^{j_i}M^i_j$, where $M^i_j$ is an antichain and for any $x\in M^i_j,y\in M^i_{j+1}$ we have $x\leqslant y$. Observe that this implies that if $x\in M^{i_1}_{j_1}\cap M^{i_2}_{j_2}$, then $M^{i_1}_{j_1+1}=M^{i_2}_{j_2+1}$. Consequently, if there exist $x_1\in M^{i_1}_{j_1}\cap M^{i_2}_{j_2}$ and $x_2 \in M^{i_2}_{j_2}\cap M^{i_3}_{j_3}$, then $M^{i_1}_{j_1+1}=M^{i_3}_{j_3+1}$, and so on. So one can merge two $M^i_j$s if they share a common element and continue this process. Then we can color $P$ as follows: all $x_i$s are Red and all other elements receive colors according to their merged class $M$. This is a proper coloring as if $x\leqslant y$ then for some $i$ and $j<j'$ we have $x\in M^i_j, y\in M^i_{j'}$ and those antichains can never be merged as $M^i_{j+1}\neq M^i_{j'+1}$. Also, this coloring does not admit a rainbow copy of $\vee$. Indeed, as only the $x_i$s are red, we can assume that the minimal element of the rainbow $\vee$ is an $x_i$. Then the two upper elements should come from different classes, but as they both belong $U_P(x_i)$, they are comparable, so they do not form a copy of $\vee$. This shows that $J$ is the only poset in $M(\vee)$.

The fact that $O_k\in F(A_k)$ follows from the blow-up algorithm of Lemma \ref{perp} that outputs $O_k$ for input $A_k$ no matter what permutation $\pi$ of $A_k$ we consider. To see that $O_k\in M(A_k)$, observe that by  Theorem \ref{mP}, we have $m(A_k)=\binom{k+1}{2}=|O_k|$, so no proper subposet of $O_k$ can rainbow force $A_k$.
\end{proof}

\begin{rem}
    The proof showing $M(A_2)=\{O_2\}$ remains valid for graphs showing that the empty graph $E_2$ on two vertices we have $M(E_2)=\{K_1+K_2\}$, where $G+F$ is the vertex-disjoint union of $G$ and $F$.
\end{rem}

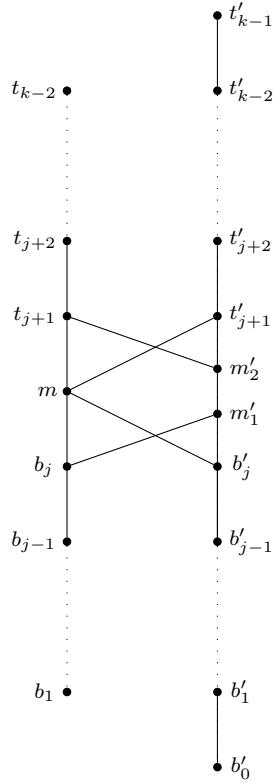
\begin{figure}
    \centering
\begin{tikzpicture}[scale=1, every node/.style={draw, circle, fill=black, inner sep=1pt}]

    \node[label=left:\scriptsize$m$] (A) at (-1,0) {};

 \draw (A);
        
        \node[label=left:\scriptsize$t_{j+1}$] (U) at (-1,1) {};
        
        \draw (U);
        \draw (A) -- (U);

        \node[label=left:\scriptsize$t_{j+2}$] (U2) at (-1,2) {};
        
        \draw (U2);
        \draw (U2) -- (U);

        \node[label=left:\scriptsize$t_{k-2}$] (Uk) at (-1,4) {};
        
        \draw (Uk);

         \draw[loosely dotted] (Uk) -- (U2);

        \node[label=right:\scriptsize$t'_{k-2}$] (U'k) at (1,4) {};
        
        \draw (U'k);

        \node[label=right:\scriptsize$t'_{k-1}$] (U'k-1) at (1,5) {};
        
        \draw (U'k);
        \draw (U'k) -- (U'k-1);
        \node (C) at (1,1) [label=right:\scriptsize$t'_{j+1}$] {};
        \draw (C);
        \draw (A) -- (C);
        \node (D) at (1,-1) [label=right:\scriptsize$b'_j$] {};
        \draw (D);
        \draw (A) -- (D);

        \node (B) at (-1,-1) [label=left:\scriptsize$b_j$] {};
        \draw (B);

        \draw (A) -- (B);
    
        \node (E) at (1,2) [label=right:\scriptsize$t'_{j+2}$] {};
        \draw (E);
        \draw (C) -- (E);
        \node (F) at (1,-2) [label=right:\scriptsize$b'_{j-1}$] {};
        \draw (F);
        \draw (D) -- (F);

        \draw[loosely dotted] (U'k) -- (F);

        \node (B2) at (-1,-2) [label=left:\scriptsize$b_{j-1}$] {};
        \draw (B2);
        \draw (B) -- (B2);

        \node[label=right:\scriptsize$m'_1$] (M') at (1,-0.3) {};
        
        \draw (M');

        \draw (B) -- (M');
        \node[label=right:\scriptsize$m'_2$] (M'') at (1,0.3) {};
        
        \draw (M'');

        \draw (F) -- (M') -- (M'') -- (C);

        \draw (M'') -- (U);

        \node[label=left:\scriptsize$b_1$] (B1) at (-1,-4) {};

        \draw[loosely dotted] (B2) -- (B1);

        \node[label=right:\scriptsize$b'_1$] (B'1) at (1,-4) {};

        \node[label=right:\scriptsize$b'_0$] (B'0) at (1,-5) {};

        \draw (B'0) -- (B'1);

         \draw[loosely dotted] (F) -- (B'0);
    \end{tikzpicture}
    \caption{Hasse diagram of $D^j_k$.}
    \label{fig:Hasse}
\end{figure}

\begin{cons}\label{ak}
For $k\ge 4$ and $1<j<k-1$, let $O^j_k$ denote the poset $O_{k-2}\oplus D^j_k$, where $D^j_k$ contains two chains $b_1\prec b_2\prec \dots \prec b_j \prec m \prec t_{j+1}\prec \dots \prec t_{k-2}$ and $b'_0\prec b'_1\prec \dots \prec b'_j\prec m'_1\prec m'_2\prec t'_{j+1}\prec \dots \prec t'_{k-1}$ with further cover relations $b'_j\prec m\prec t'_{j+1}$ and $b_j\prec m'_1\prec m'_2\prec t_{j+1}$. (see Figure \ref{fig:Hasse})
\end{cons}

\begin{prop}\label{ok}
    For $k\ge 4$, we have $|M(A_k)|\ge k-3$.
\end{prop}

\begin{proof}
It is enough to show that $O^j_k\in M(A_k)$ for all values of $j,k$ in Construction \ref{ak}. \begin{claim} 
$O^j_k$ rainbow forces $A_k$. 
\end{claim}

\begin{proof}[Proof of Claim]
    Let $c$ be a proper coloring of $O^j_k$. As $O_{k-2}$ rainbow forces $A_{k-2}$, we find a rainbow copy of $A_{k-2}$ in the $O_{k-2}$ part of $O^j_k$. The colors used at the $A_{k-2}$ are called \textit{old colors}, while other colors are \textit{new colors}. Also, $B,B',T$, and $T'$ denote the sets of $b$s, $b'$s, $t$s, and $t'$s in the definition of $O^j_k$.

    Assume first that $T$ admits at least 2 new colors. Then as $O_2\in M(A_2)$ and $T \perp T'$, either we obtain a rainbow $A_2$ using only new colors and so a rainbow $A_k$, or $T'$ only admits old colors. But then as $B'\cup T'$ is a chain, $B'$ must admit at least 2 new colors, and thus either we find a rainbow $A_k$ or $B$ should admit only old colors. But $B\cup T'$ is a chain of length $k-1$ and there are only $k-2$ old colors, so this is impossible, and therefore we must have previously obtained a rainbow $A_k$.

    Assume next $T$ admits exactly 1 new color. Then as before, either $T'$ has at most one new color, or we obtain a rainbow $A_k$. If $B'$ admits at most one new color, then $B'\cup T'$ admit all old colors, and so $m,m'_1,m'_2$ have all new colors and we obtain a rainbow $A_k$. So $B'$ admits at least 2 new colors. But then as before, either we are done, or $B$ admits only old colors, so $B\cup T'$ has all old colors and again $m,m'_1,m'_2$ have only new colors and we obtain a rainbow $A_k$.

    Finally, suppose $T$ admits only old colors. Then if $B$ has a new color, then we can proceed as in one of the twwo previous cases since he role of $B$ and $T$ are similar. So we can assume $B\cup T$ admits only old colors. Then $B\cup T$ has all old colors, and again $m,m'_1,m'_2$ have new colors and we obtain a rainbow $A_k$.
\end{proof}

We are left to show that for any $p\in O^j_k$, the poset $O^j_k\setminus \{p\}$ does not rainbow force $A_k$. Observe that the $O_{k-2}$ part of $O^j_k$ together with $B\cup \{m\}\cup T$ induces an $O_{k-1}$ in $O^j_k$. So if $p\in O_{k-2}\cup B\cup \{m\}\cup T$, then there is a coloring without a rainbow $A_{k-1}$ and no matter how we extend this to $B'\cup \{m'_1,m'_2\}\cup T'$ properly, we can only add at most one element from $B'\cup \{m'_1,m'_2\}\cup T'$ to an antichain, so no rainbow $A_k$.

So we are left with $p\in B'\cup \{m'_1,m'_2\}\cup T'$. If $p=m'_a$, then we can consider the coloring $c$ that assigns colors $1,2,\dots,i$ to the $i$th chain of $O_{k-2}$ for all $1\le i\le k-2$, color $h$ to $b_h,b'_h,t_h,t'_h$ for $1\le h\le k-2$ and Red to $b'_0$, Blue to $t'_{k-1}$ and Green to $m,m'_{3-a}$. This does not admit any rainbow $A_k$ as for a rainbow $A_{k-2}$ in the $O_{k-2}$ part, we have to use all colors $1,2,\dots,k-2$, so we must need $m$ to finish, but it is only incomparable to $m'_{3-a}$ that received the same color as $m$.

The $p\in T'$ and $p\in B'$ cases are analogous, so we can assume $p\in T'$ and as the remaining poset is the same for all $h$ we can assume $p=t'_{k-1}$. Then for any $1\le h \le k-3$ color $b_h,b'_h,t_h,t'_h$ with $h$, the chains of $O_{k-2}$ as in the previous case, color $b'_0$ with Red, color $t_{k-2},t'_{k-2}$ Blue, color $m,m'_1$ Green, and finally color $m'_2$ with $k-2$. It is not hard to see that this coloring does not admit a rainbow $A_k$.
\end{proof}

\begin{rem}
    It is not hard to show that for any $4\le k \le \ell$ and $1<j<k-1$, the poset $O^j_k\oplus \bigoplus_{i=k+1}^\ell C_i\in M(A_\ell)$. More generally, if in $O_{2k}$ for some $I\subseteq \{2,3,\dots,k\}$, we replace $C_{2i-1}\oplus C_{2i}$ by $D^{j_i}_{2i}$ for all $i\in I$, then we obtain a poset in $M(A_{2k})$. Also, $|M(A_3)|>1$ can be shown by observing that the poset on 7 elements $u,v_1,v_2,w_0,w_1,w_2,w_3$ with $v_1\prec v_2, v_1\prec w_2, w_1\prec v_2, w_0\prec w_1\prec w_2 \prec w_3$ being all its covering relations minimally rainbow forces $A_3$.
\end{rem}

\begin{rem}
For the empty graph $E_k$ on $k$ vertices, we have that all long enough odd cycles belong to $M(E_3)$, and for any $k>3$ a long enough odd cycle together with $\cup_{i=4}^kK_k$ are graphs in $M(E_k)$. Therefore, $M(E_k)$ is infinite for all $k\ge 3$.
\end{rem}

\begin{prob}
Decide whether $M(A_3)$ is finite or infinite.    
\end{prob}

\section{Results on $\RLA(n,P)$}

In this section, we prove results on $\RLA(n,P)$ for various posets $P$. We start with the proof of Proposition \ref{connection} that establishes the connection between the functions $\LA$ and $\RLA$.

\begin{proof}[Proof of Proposition \ref{connection}]
    Assume first $\cF\subseteq 2^{[n]}$ with $|\cF|>\LA(n,M(P))$. Then $\cF$ contains a copy of some $Q\in M(P)$, and so as $Q$ rainbow forces $P$, therefore $\cF$ cannot be colored properly without admitting a rainbow copy of $P$. This shows $\RLA(n,P)\le \LA(n,M(P))$.

    Assume next $\cF\subseteq 2^{[n]}$ with $|\cF|>\RLA(n,P)$. Then no matter how we properly color $\cF$, the coloring admits a rainbow copy of $P$. So $(\cF,\subseteq)$ rainbow forces $P$ and therefore $(\cF,\subseteq)$ contains a member $Q$ of $M(P)$ and so $\LA(n,M(P))\le \RLA(n,P)$. 
\end{proof}

\subsection{Tree posets} In this subsection, we prove Theorem \ref{tree}. We need some definitions and notations.

A tree poset is a \textit{downtree} if it has a largest element, and an \textit{uptree} if it has a smallest element. Let $\Delta^{k,h}$ denote the the downtree with $\overrightarrow{H}(\Delta^{k,h})$ being the complete $k$-ary tree of height $h$ (here height means the number of vertices in a longest path to the root) with all edges pointing towards the root, and let $\nabla^{k,h}$ denote the dual poset of $\Delta^{k,h}$. We write $\Delta^k=\Delta^{k,k}$ and $\nabla^k=\nabla^{k,k}$. 
For a coloring $c$ of $P$, $x\in P$, and a set $S$ of colors, we write $S_c(x)=\{u\in S: \exists y \in D_P(x) ~c(y)=u\}$ and $S^c(x)=\{u\in S: \exists y \in U_P(x) ~c(y)=u\}$. The next lemma is crucial in the proof of Theorem \ref{tree}. Let us remind the reader that the rank of an element $x$ in a poset is the length of the longest chain ending in $x$. So the root of $\Delta^k$ has rank $k$, its children in the Hasse diagram have rank $k-1$, and so on.

For a downtree or uptree $T$, we say that $f:T\rightarrow |T|$ is a \textit{finishing order} if there exists a depth-first-search of $T$ starting from its largest or smallest element for which the finishing order is $f$. Recall that if $y$ is the root of a downtree $T$ and $y_1,y_2\dots,y_s$ are the children of $y$ in $T$ with $f(y_1)<f(y_2)<\dots<f(y_S)$, then the finishing numbers of elements in $D_T(y_j)$ form the interval $[f(y_{j-1})+1,f(y_j)]$. 

In the next lemma and theorem, we are going to define an embedding of $T$ into $\Delta^k$ and to the to-be-defined $\cT^k$. For any $t\in T$, $z_t$ will denote its image.

\begin{lem}\label{downlemma}
    For any downtree $T$, finishing order $f$, a set $K$ of (unusable) colors, and $k\ge |T|+|K|$ the following holds: in any proper coloring $c$ of $\Delta^k$ there exists a rainbow copy of $T$ such that
    \begin{enumerate}
        \item 
        the largest element of $T$ is placed in a highest ranked element $x$ of $\Delta^k$ with $c(x)\notin K$,
        \item 
        for any element $z$ of the rainbow copy of $T$, we have $c(z)\notin K$,
        \item 
        for any non-root element $t\in T$, the corresponding element $z_t$ in the rainbow copy satisfies $r(z_t)\le |K(f,t)_c(z_t)|+1$, where $K(f,t)=K\cup \{c(z_{t'}):f(t')<f(t)\}$.
    \end{enumerate}
\end{lem}

\begin{proof}
    We proceed by induction on $|T|+|K|$ with the base case $|T|+|K|=1$ being trivial as it means $K$ is empty, $T$ is the one-element poset, and thus one can place this element to the root of $\Delta^k$ in any coloring $c$.

    So fix $T$ and $K$, and suppose $c$ is a proper coloring of $\Delta^k$. We can assume that the root $x_r$ has a color not in $K$, as otherwise $c(x_r)$ does not appear on any other vertex of $\Delta^k$, so one can consider a subposet of $D_{\Delta^k}(x)$ isomorphic to $\Delta^{k-1}$ for any child $x$ of $x_r$ with $K'=K\setminus \{c(x_r)\}$. $k\ge |T|+|K|$ implies $k-1\ge |T|+|K'|$ so one can apply induction to $D_{\Delta^k}(x)$ and $K'$.

    Therefore, we can assume that we can place the largest element $y$ of $T$ into the root $x$ of $\Delta^k$. As $c$ is a proper coloring, $c(x)$ does not appear on any other element of $\Delta^k$. Let $y_1,y_2,\dots,y_j$ be the children of $y$ in $T$ with $f(y_1)<f(y_2)<\dots<f(y_s)$. We plan to use induction and find rainbow copies of $D_T(y_i)$ in $D_{\Delta^k}(x_i)$, where $x_i$ is the $i$th child of $x$. (Note that $j\le |T|-1\le k$, so $x$ has enough children.) Assume that we have managed to find rainbow copies of $D_T(y_i)$ for $i=1,2,\dots,\ell$ such that the colors used are all distinct (so their union together with the root is also rainbow) and satisfy  properties (2) and (3) of the lemma. Then let $K^\ell$ be the color set of $K$ together will all colors used in the copies of $D_T(y_1), D_T(y_2), \dots, D_T(y_\ell)$. As those copies are rainbow, we have $|K^\ell|=|K|+\sum_{i=1}^\ell|D_T(y_i)|=|K|+f(y_\ell)$. As $\sum_{i=1}^j|D_T(y_i)|=|T|-1$, we obtain that $|K^\ell|+|D_T(y_{\ell+1})|=|K|+\sum_{i=1}^{\ell+1}|D_T(y_i)|\le |K|+|T|-1\le k-1$, so we could apply induction to $K^\ell$, $D_T(y_{\ell+1})$, the finishing order $f'(t):=f(t)-f(y_\ell)$, and a copy of $\Delta^{k-1}$ in $D_{\Delta^k}(x_{\ell+1})$. We obtain a rainbow copy of $D_T(y_{\ell+1})$ such that for any $t\in D_T(y_{\ell+1})\setminus \{y_{\ell+1}\}$ and for its image $z_t$ we have $r(z_t)\le |K^{\ell,f',t}_c(z_t)|+1\le |K^{f,t}_c(z_t)|+1$. The last inequality follows from $K\subset K^\ell$ and $|K^\ell|-|K|= f(y_\ell)$, $f'(t)=f(t)-f(y_\ell)$. So the only problem is that because of property (1), the image of $y_{\ell+1}$ would be placed possibly to $x_{\ell+1}$ and thus might not satisfy property (3).

    Observe that for any $z\in \Delta^k$ and its child $z'$, we have $r(z')=r(z)-1$ and $|K^\ell_c(z)|\ge |K^\ell_c(z')|$. Therefore, $r(z')-1-(|D_T(y_{\ell+1})|+|K^\ell_c(z')|)$ can be at most one smaller than $r(z)-1-(|D_T(y_{\ell+1})|+|K^\ell_c(z)|)$ (and the former can also be equal to or larger than the latter). Also, if $w$ is a leaf of $\Delta^k$, then $r(w)-1-(|D_T(y_{\ell+1})|+|K^\ell_c(w)|)$ is negative. So either $r(x_{\ell+1})-1\le |D_T(y_{\ell+1})|+|K^\ell_c(x_{\ell+1})|$ or, by applying the above observations repeatedly to consecutive elements of a chain from $x_{\ell+1}$ to a leaf, there exists $x'_{\ell+1}\in D_{\Delta^k}(x_{\ell+1})$ with $r(x'_{\ell+1})-1=|D_T(y_{\ell+1})|+|K^\ell_c(x'_{\ell+1})|$. In either case, we can apply induction to $D_{\Delta^k}(x_{\ell+1}), f'$ or to $D_{\Delta^k}(x'_{\ell+1}),f'$ and because of the extra condition on $x_{\ell+1}$ or $x'_{\ell+1}$, property (3) of the lemma will be satisfied for the image of $y_{\ell+1}$, too.
\end{proof}

By changing to the dual in Lemma \ref{downlemma}, one obtains the following statement.

\begin{lem}\label{uplemma}
    For any uptree $T$, finishing order $f$, a set $K$ of (unusable) colors, and $k\ge |T|+|K|$ the following holds: in any proper coloring $c$ of $\nabla^k$ there exists a rainbow copy of $T$ such that
    \begin{enumerate}
        \item 
        the smallest element of $T$ is placed in a lowest ranked element $x$ of $\Delta^k$ with $c(x)\notin K$,
        \item 
        for any element $z$ of the rainbow copy of $T$, we have $c(z)\notin K$,
        \item 
        for any non-root element $t\in T$, the corresponding element $z_t$ in the rainbow copy satisfies $k-r(z_t)+1\le |K(f,t)^c(z_t)|+1$, where $K(f,t)=K\cup \{c(z_{t'}):f(t')<f(t)\}$.
    \end{enumerate}
\end{lem}

With Lemma \ref{downlemma} and \ref{uplemma} in hand, we are going to define a tree poset $\cT^k$ that will be universal for all tree posets of size $k$ in the sense that $\cT^k$ will rainbow force all such tree posets. The poset $\cT^k$ has height $k$ and all its maximal chains have length $k$, $CG(\cT^k)$ has radius $k$ with its central vertex having rank $k$. We generate $\cT^k$ is as follows. Consider a copy of $\Delta^k$, its root $x$ is going to be the central element of $\cT^k$. The other elements of this $\Delta^k$ are \textit{down elements}. To any such element $y$ with $r(y)=r$, we append a $\nabla^{k,k-r+1}$ identifying $y$ with the root of this uptree and denote this by $U'(y)$. All other elements of these uptrees are \textit{up elements}. Note that these elements are at distance 2 from $x$ in $CG(\cT^k)$. We continue in this way: to every up element $z$ of rank $r$, we append a $\Delta^{k,r}$, denoted $D'(z)$, and make the newly added elements down elements, while to every down element $z$ of rank $r$, we append a $\nabla^{k,k-r+1}$, denoted $U'(z)$, and make the newly added elements up elements. We finish this process when the radius of $CG(\cT^k)$ becomes $k$, i.e. after $k$ steps.

\medskip

We need some preliminary notation for $T$ in order to present the proof of the next theorem that, as mentioned in the Introduction, will immediately imply Theorem \ref{tree}. Let us fix a maximal element $y\in T$. Elements that have odd distance from $y$ in $CG(T)$ are \textit{down elements} of $T$, while elements of even distance from $y$ in $CG(T)$ are \textit{up elements} of $T$. For any $t\in T$, let us write $\dist(t)$ to denote the distance of $t$ from $y$ in $CG(T)$. For an up element $t\in T$, let its \textit{originator} be the element $t'$ that maximizes $r(t')$ over all elements in $D_T(t)$ with $\dist(t')=\dist(t)-1$. Similarly, the originator of a down element $t$ is the element $t'$ that minimizes $r(t')$ over all elements in $U_T(t)$ with $\dist(t')=\dist(t)-1$. We denote the originator of $t$ by $\orig(t)$. Finally, let us define $U'_T(t)=\{t':\orig(t')=t\}\cup \{t\}$ for any down element $t$, and $D'_T(t)=\{t':\orig(t')=t\}\cup \{t\}$ for any up element $t$. Observe that $U'_T(t)\subseteq U_T(t)$ for down elements and $D'_T(t)\subseteq D_T(t)$ for up elements. Also, \[
T=D_T(y)\cup \bigcup_{t~\text{is up}}(D'_T(t)\setminus \{t\})\cup \bigcup_{t~\text{is down}}(U'_T(t)\setminus \{t\})
\] is a partition of $T$.

\begin{thm}\label{auxtree}
    For any tree poset $T$ on $k$ elements, $\cT^k$ rainbow forces $T$.
\end{thm}

\begin{proof}
    Fix an enumeration $t_1,t_2,\dots,t_{|T|}$ such that $i\le j$ implies $\dist(t_i)\le \dist(t_i)$, in particular, $t_1=y$. 
    Let $c$ be a proper coloring of $\cT^k$. We will find a rainbow copy of $T$ in $\cT^k$ by applying Lemma \ref{downlemma} or Lemma \ref{uplemma} to $D'_T(t_j)$ or $U'_T(t_j)$ according to the fixed ordering of the elements of $T$. We first let $z_y=z_{t_1}=x$. Then we can apply Lemma \ref{downlemma} to $D_T(y)$, any finishing ordering $f$ of $D_T(y)$ and $K^1=\emptyset$.

    Assume that for all $i\le j$, we were able to embed $U'_T(t_i)$ for all down elements and $D'_T(t_i)$ for all up elements such that the embedding used Lemma \ref{uplemma} for down elements, and Lemma \ref{downlemma} for up elements, and for any element $t_i$, the unusable color set $K^i$ was a subset of \[L_i=\left\{c(z_{t_h}): t_h\in D_T(y)\cup \bigcup_{t_\ell ~\text{is up,}\ \ell<i}D'_T(t_\ell) \cup \bigcup_{t\ell ~\text{is down,}\ \ell<i}U'_T(t_\ell)\right\}.
    \]
    That is, $K^i$ was a subset of the colors used at elements of the already defined part of the copy of $T$ at that moment.

    We claim that we can apply Lemma \ref{downlemma} or Lemma \ref{uplemma} (depending on whether $t_j$ is an up or a down element) to embed $D'_T(t_j)$ or $U'_T(t_j)$. Let $t_m=\orig(t_j)$ for some $m<j$. Assume $t_m$ is an up element (or $y$) and $t_j$ is a down element, the other case is analogous. There are three kinds of colors in $L_j$:
    \begin{enumerate}
        \item 
       Colors that appear below $z_{t_j}$.
        \item 
        Colors that do not appear below $z_{t_j}$.
        \item 
        The color of $z_{t_j}$.
    \end{enumerate}
    We let $K^j$ be the set of colors of type (2). Lemma \ref{downlemma} implies that $r(t_j)-1$ is at most the number of colors appearing below $z_{t_j}$ that are also colors of elements embedded no later than $D'_T(t_m)$. So $r(z_{t_j})-1$ is a most the number of colors of type (1). Colors of type (1) need not be included in $K^j$ as $c$ is a proper coloring, these colors cannot appear above $z_{t_j}$. By definition, we have $|L_i|-1=\#(\text{type (1)})+\#(\text{type (2)})$. Recall that $U_{\cT^k}(z_{t_j})=\nabla^{k,k-(r(z_{t_j})-1)}$. We want to embed $U'_T(t_j)$ that has size at most 
    \[
    |T|-|L_i|+1=|T|-(\#(\text{type (1)})+\#(\text{type (2)}))=|T|-(\#(\text{type (1)}+|K^j|).
    \]
    So 
    \[
    |U'_T(t_j)|+|K^j|\le |T|-\#(\text{type (1)})\le |T|-(r(z_{t_j})-1),
    \]
    the height of $U_{\cT^k}(z_{t_j})$. So the conditions of Lemma \ref{uplemma} are indeed satisfied for a subposet $\nabla^{k-(r(z_{t_j})+1}$ of $U_{\cT^k}(z_{t_j})$, any of its finishing order $f$ and $K=K^j$. This finishes the proof of the theorem.
\end{proof}

\begin{proof}[Proof of Theorem \ref{tree}]
    For the lower bound, one can consider the middle $|T|-1$ layers colored by layers. As there is not enough colors, there is no rainbow copy of $T$.

    For the upper, by Theorem \ref{auxtree} and Theorem \ref{bj}, any family $\cF\subseteq 2^{[n]}$ of size $(|T|-1+\varepsilon)\binom{n}{\lfloor n/2\rfloor}$ contains a copy of $\cT^{|T|}$ which rainbow forces $T$, so any proper coloring of $\cF$ admits a rainbow copy of $T$.
\end{proof}

\subsection{Multi-partite posets and antichains} In this subsection, we prove Theorem \ref{diamond} and Proposition \ref{easylar}. We start by defining the notions that we will use in our proofs.

Let $\mathbf{C}_n$ denote the set of all $n!$ maximal chains in $[n]$. The \textit{Lubell-mass} of a family $\cF\subseteq 2^{[n]}$ is defined as
\[
\lambda_n(\cF)=\sum_{F\in \cF}\frac{1}{\binom{n}{|F|}}=\frac{1}{n!}\sum_{\cC\in \mathbf{C}_n}|\cC\cap \cF|,
\]
the average number of sets in $\cF$ that a maximal chain contains when $C$ is chosen uniformly at random from $\mathbf{C}_n$. The famous LYM-inequality states that the Lubell-mass of an antichain is at most 1.

\begin{thm}[LYM-inequality \cite{Lub,Mesh,Yam}]\label{lym}
    For any antichain $\cA\subseteq 2^{[n]}$, we have $\lambda_n(\cA)\le 1$ with equality only for $\cA=\binom{[n]}{i}$ for some $i=0,1,\dots,n$.
\end{thm}

Mirsky's theorem \cite{M} states that if the largest chain in a poset $P$ contains $k$ elements, then $P$ can be partitioned into $k$ antichains. This immediately implies the following.

\begin{cor}\label{klym}
    If $\cF\subseteq 2^{[n]}$ does not contain chains of length $k+1$, then $\lambda_n(\cF)\le k$.
\end{cor}

Families satisfying the condition of Corollary \ref{klym} are called \textit{$k$-Sperner} (and antichains are also called Sperner families). An immediate consequence of Corollary \ref{klym} is a theorem of Erd\H os \cite{E} stating that any $k$-Sperner family $\cF\subseteq 2^{[n]}$ has size at most $\Sigma(n,k)$. We will need the following very weak stability result.

\begin{cor}\label{stab}
    If a $k$-Sperner family $\cF\subseteq 2^{[n]}$ contains a set of size $i<n/2$, then $|\cF|\le \Sigma(n,k)-\frac{\binom{n}{\lfloor \frac{n-k}{2}\rfloor}}{\binom{n}{i}}+1$.
\end{cor}

\begin{proof}
    The size of $\cF$ is the number of summands in $\lambda_n(\cF)$. The smallest summands correspond to sets close to size $n/2$. A set of size $i$ yields a summand $\frac{1}{\binom{n}{i}}$ which is $\frac{\binom{n}{\lfloor \frac{n-k}{2}\rfloor}}{\binom{n}{i}}$ times larger than the largest summand in the extremal family.
\end{proof}

\begin{lem}[Griggs, Li, Lu, Lemma 3.2 in \cite{GLL}]\label{lubm}
    For any family $\cF\subseteq 2^{[n]}$, we have $|\cF|\le \lambda_n(\cF)\binom{n}{\lfloor n/2\rfloor}$.
\end{lem}

\begin{lem}\label{h2lym}
    If $\cF\subseteq 2^{[n]}$ is a $H_2$-free family with $\emptyset,[n]\in \cF$, then $\lambda_n(\cF)\le 3+\frac{1}{6}$ and the bound is sharp.
\end{lem}

\begin{proof}
    Observe first that taking $n=4$, for the family $\cF_4=\{\{\emptyset\},\{1\},\{2\},\{1,2\},\{1,2,3\},\{1,2,4\},[4]\}$ we have $\lambda_4(\cF_4)=3+\frac{1}{6}$. Note that $\cF_4$ is the family of all subsets of $[4]$ that are comparable to $\{1,2\}$.

    Consider next an $H_2$-free family $\cF\subseteq 2^{[n]}$ with $\emptyset,[n]\in \cF$. Then $\cF':=\cF\setminus\{\emptyset,[n]\}$ is $O_2$-free and we need to show $\lambda_n(\cF')\le 1+\frac{1}{6}$. As in the proof of Lemma \ref{smallexample}, the poset structure of $\cF'$ must be complete multipartite $K_{s_1,s_2,\dots,s_k}$. Moreover, if $s_i,s_{i+1}>1$, then one can add the union of the sets corresponding to level $i$ and the intersection of the sets of level $i+1$ (these two might coincide). We proceed by induction on $n+k$. If $k=1$, then $\cF'$ is an antichain, and by Theorem \ref{lym}, we have $\lambda_n(\cF')\le 1$. If $n=2$, then $\lambda_2(2^{[2]})=3$ or if $n=3$, then $\cF'$ can contain at most 3 sets, as if it contains at least two singletons, then it might further contain either their union or the third singleton, but not both. So again $\lambda_3(\cF')\le 1$. 

The inductive step: we consider two cases according to $s_k=1$ or $s_k>1$. Assume first $s_k=1$ with the largest set $M$ of $\cF'$ having size $m\le n-1$. Then by induction, we obtain $$\lambda_m(\cF'\setminus\{M \})=\sum_{F \in \cF'\setminus \{M\}}\frac{1}{\binom{m}{|F|}}\le 1+\frac{1}{6}.
$$
As for any $1\le i\le m-1$, we have $\frac{\binom{m}{i}}{\binom{n}{i}}\le \frac{m}{n}$, we reach
\[
\lambda_n(\cF')=\sum_{F\in \cF'}\frac{1}{\binom{n}{|F|}}\le \frac{1}{\binom{n}{|M|}}+\frac{m}{n}\sum_{F \in \cF'\setminus \{M\}}\frac{1}{\binom{m}{|F|}}\le \frac{1}{n}+\left(1-\frac{1}{n}\right)\left(1+\frac{1}{6}\right)\le 1+\frac{1}{6}
\]
as wanted.

Assume next, $s_k>1$ and thus by the above observation, $s_{k-1}=1$. If $k=2$, then we can consider $\overline{\cF'}=\{[n]\setminus F: F\in \cF'\}$. As $\lambda_n(\cF')=\lambda_n(\overline{\cF'})$, by the previous case, we have $\lambda_n(\cF')\le 1+\frac{1}{6}$. So we can assume $k\ge 3$. Let $m$ be the size of $M$, the set in $\cF'$ corresponding to its $(k-1)$st part, and note that $2\le m \le n-2$. The complements of the sets in $\cF'$ properly containing $M$ form an antichain in $[n]\setminus M$, so we obtain $\sum_{F\in \cF':M\subset F}\frac{1}{\binom{n-m}{n-|F|}}\le 1$. Therefore, applying induction
\[
\lambda_n(\cF')=\frac{1}{\binom{n}{m}}+\sum_{F\in \cF':F\subset M}\frac{1}{\binom{n}{|F|}}+\sum_{F\in \cF':M\subset F}\frac{1}{\binom{n}{n-|F|}}\le \frac{2}{n(n-1)}+\frac{m}{n}\cdot\frac{7}{6}+\frac{n-m}{n}\le \frac{7}{6}-\frac{1}{6}\cdot\frac{2}{n}+\frac{2}{n(n-1)},
\]
which is at most $1+\frac{1}{6}$ if $n$ is at least 7.

If $k=3$, then the sets in $\cF'$ contained in $M$ together with the complements of the sets in $\cF'$ that contain $M$ form an antichain in $2^{[n]}$, so by Theorem \ref{lym}, we have $\lambda_n(\cF')\le 1+\frac{1}{\binom{n}{m}}\le 1+\frac{1}{\binom{n}{2}}\le \frac{7}{6}$ with equality only if $m=2$ and $n=4$, as in the construction showing tightness. So we can assume $k\ge 4$. If $k=4$ and $n=5$, then the poset structure of $\cF'$ is either $K_{2,1,1,2}$ or $K_{1,2,1,2}$. In the former case, $\lambda_n(\cF')=\frac{4}{5}+\frac{2}{10}$, while in the latter case $\lambda_n(\cF')=\frac{3}{5}+\frac{3}{10}$, both at most 1. Finally, if $n=6$, then $k=4$ and $k=5$ are possible. If $k=5$, then the poset structure $K_{2,1,2,1,2}$ gives largest Lubell-mass $\frac{4}{6}+\frac{2}{15}+\frac{2}{20}$. If $k=4$, then the poset structure can be $K_{a,1,b,1}$, $K_{1,a,1,b}$, or $K_{a,1,1,b}$. For the sets $F\subsetneq F'$ we have $1\le |F|<|F'|\le 5$ and by their positioning, we have that among $|F|, |F'|-|F|,n-|F'|$ the largest can be 3 and all others can be at most 2. As the largest antichain in $2^{[3]}$ and $2^{[2]}$ have size 3 and 2 respectively, we have that $a+b\le 2+3$. So there are at most most 5 sets of size 1 or 5 and a total of 7 sets giving Lubell-mass at most $\frac{5}{6}+\frac{2}{15}$. This finishes the proof of the lemma.
\end{proof}

Theorem \ref{diamond} will now easily follow from Lemma \ref{lubm}, Lemma \ref{h2lym}, and the chain-partition method introduced by Griggs and Li \cite{GL}.

\begin{proof}
    Let $\cF\subseteq 2^{[n]}$ be a rainbow $\Diamond$-free family. Then by Lemma \ref{smallexample} and Proposition \ref{connection}, $\cF$ is $H_2$-free. Let us introduce the min-max partition of $\textbf{C}_n$: for any pair $F\subsetneq F'$ with $F,F'\in \cF$, we denote by $\textbf{C}_{F,F'}$ the collection of those maximal chain $\cC\in \textbf{C}_n$ for which the smallest set in $\cF\cap \cC$ is $F$ and the largest set in $\cF\cap \cC$ is $F'$. Let us write $\textbf{C}^-:=\textbf{C}_n\setminus \cup_{F\subsetneq F'}\textbf{C}_{F,F'}$. By definition, for any $\cC\in \textbf{C}^-$, we have $|\cC\cap \cF|\le 1$. Also by definition, for any $F\subseteq F'$ and maximal chain $\cC\in \textbf{C}_{F,F'}$, $\cC$ does not meet $\cF$ 'below' $F$ nor 'above' $F'$. By Lemma \ref{h2lym}, the average size of $\cC\cap \cF$ between $F$ and $F'$ is at most $3+\frac{1}{6}$. As $\lambda_n(\cF)$ is a weighted average of the average size of $\cC\cap \cF$ in $\textbf{C}_{F,F'}$ and in $\textbf{C}^-$, we obtain that $\lambda_n(\cF)\le 3+\frac{1}{6}$. Lemma \ref{lubm} implies the statement of the theorem. 
\end{proof}

\medskip

Finally, we prove our results on $\RLA(n,A_k)$ and $\RLA(n,K_{s,t})$.

\begin{proof}[Proof of Proposition \ref{easylar}]
    The lower bound of part (i) follows from considering the family of the middle $k-1$ layers of $2^{[n]}$ with $\emptyset$ and $[n]$ added. As $\emptyset$ and $[n]$ are related to all sets, they cannot be part of any $A_k$. The rest of the family can be properly colored by layers using only $k-1$ colors, so without rainbow $A_k$. This shows $\RLA(n,A_k)\ge \Sigma(n,k-1)+2$.

    To see the upper bound of (i), assume $\cF\subseteq 2^{[n]}$ is of size $\Sigma(n, k-1)+3$. We will show by induction on $k$ that $\cF$ contains an $O_k$, so by Proposition \ref{smallexample}, $\cF$ contains a rainbow $A_k$ for any proper coloring. As  base case, we will use $k=1$, as then $1=\Sigma(n,1-1)+1\le \Sigma(n,1-1)+3$ set is enough for a rainbow $A_1$.
    
    Let us write $i_n=\lfloor n/10\rfloor$. This choice implies that for large enough $n$, we have 
    \begin{equation}\label{hin}
        4\binom{n}{i_n}^2\le 2^{2h(1/10)n}<\binom{n}{\lfloor \frac{n-k}{2}\rfloor},
    \end{equation} 
    where $h$ is binary entropy function $h(c)=-\log_2c-(1-c)\log_2(1-c)$ and $h(1/10)<0.47$. Furthermore
    \begin{equation}\label{geometric}
        \sum_{j=0}^{i_n}\binom{n}{j}\le \left(\frac{9}{8}+o(1)\right)\binom{n}{i_n}. 
    \end{equation}
    Let us partition $\cF\setminus \{\emptyset,[n]\}$ into $\cF^{mid}\cup \cF'$ with $\cF^{mid}=\{F\in \cF: i_n<|F|<n-i_n\}$ and $\cF'=\cF\setminus (\{\emptyset,[n]\}\cup \cF^{mid})$. Note that if $n$ is large enough, then $|\cF'|=O(2^{h(1/10)n})\le \frac{1}{n}\Sigma(n,k-1)$. Suppose first that $\cF^{mid}$ contains a $C_k$. Then the total number of sets in $2^{[n]}$ that are comparable to any set in this $C_k$ is at most $2\cdot 2^{0.9n}$. Let us denote the family of these sets by $\cG$. As $|\cF\setminus\cG|\ge \Sigma(n,k-1)+3-2\cdot 2^{0.9n}>\Sigma(n,k-2)+3$ if $n$ is large enough, therefore, by induction, $\cF\setminus \cG$ contains a copy of $O_{k-1}$, and so together with $C_k$, $\cF$ contains a copy of $O_k$ as claimed.

    So we can assume that $\cF^{mid}$ is $(k-1)$-Sperner and thus by Erd\H os's theorem, we have $|\cF^{mid}|\le \Sigma(n,k-1)$. If $\cF'$ is empty, then $|\cF|\le \Sigma(n,k-1)+2$, so $\cF'$ is non-empty. Suppose next that there exists $F\in \cF'$ such that $\{F\}\cup \cF^{mid}$ is also $(k-1)$-Sperner. Then by (\ref{geometric}), we have $|\cF'|\le 2\sum_{j=0}^{i_n}\binom{n}{j}=(\frac{9}{4}+o(1))\binom{n}{i_n}$. Applying this and Corollary \ref{stab} with (\ref{hin}), we obtain
    \[
    |\cF|=|\{\emptyset,[n]\}|+|\cF^{mid}|+|\cF'|\le 2+\Sigma(n,k-1)-4\binom{n}{i_n}+\left(\frac{9}{4}+o(1)\right)\binom{n}{i_n}<\Sigma(n,k-1)
    \]
    which is a contradiction.
    
    Finally, we consider the case that for some $F\in \cF'$, $\{F\}\cup \cF^{mid}$ contains a $C_k$. By taking complements if necessary, we can assume that $F$ is the bottom element of $C_k$. Let $\cG$ denote again the family of sets that are comparable to at least one set in $C_k$. Fix $z\in F$. We want a lower bound on $|\{H\subseteq [n]:z\notin H\}\cap \cF^{mid}\setminus \cG|$. The number of sets that are below the top element of $C_k$ is at most $2^{0.9n}$ as the top element belongs to $\cF^{mid}$. As $\cF^{mid}$ is $(k-1)$-Sperner, the number of sets in $\cF^{mid}$ that contain $z$ is at most $\Sigma(n-1,k-1)$. As $z\in F$, all sets containing $F$ contain $z$. Therefore, using $\Sigma(n,k-1)=(2+o(1))\Sigma(n-1,k-1)$ and $\Sigma(n,k-1)=\Theta(\frac{1}{\sqrt{n}}2^n)$, we obtain 
    \begin{align*}
        |\{H\subseteq [n]:z\notin H\}\cap \cF^{mid}\setminus \cG| & \ge \Sigma(n,k-1)-2^{0.9n}-\Sigma(n-1,k-1)-\frac{1}{n}\Sigma(n,k-1)  \\
         & =(1-o(1))\Sigma(n-1,k-1)\ge \Sigma(n-1,k-2)+3.
    \end{align*} 
    Thus by induction $\cF^{mid}\setminus \cG$ contains an $O_{k-1}$ and $\cF$ contains an $O_k$.

    \medskip
    
    The lower bound of part (ii) follows from considering the family of the middle $s+t$ layers of $2^{[n]}$ properly colored according to set sizes. We claim it does not contain rainbow copies of $K_{s,t}$. Indeed, by the rainbow property, all sets representing an element of $K_{s,t}$ should come from distinct layers. The intersection of the sets representing the top elements is smaller than all of these representatives and larger than the representatives of all the $s$ bottom elements. That would require one more layer. 

    To see the upper bound of (ii), observe that as established in the previous section, $O_k\in M(A_k)$. Then by Lemma \ref{linsum}, $O_s\olessthan O_t\in M(A_s\olessthan A_t)$ and $A_s\olessthan A_t=K_{s,t}$. But $O_s\olessthan O_t$ is a subposet of $(O_s\olessthan C_1)\olessthan O_t$, which in turn is a tree-poset of height $s+t+1$. So by Proposition \ref{connection}, $\RLA(n,K_{s,t})\le \LA(n,(O_s\olessthan C_1) \olessthan O_t)=(s+t+o(1))\binom{n}{\lfloor n/2\rfloor}$, where in the last equation we used Theorem \ref{bj}.
\end{proof}

\section{Concluding remarks}

The reader who is not familiar with the area of forbidden subposet problems might wonder why our notation uses $^*$s. In the literature $\LAw(n,P)$ and $e(P)$ correspond to \textit{weak copies of $P$.} A family $\cF$ is a weak copy of $(P,\leqslant)$ if there exists a bijection $b:P\rightarrow \cF$ such that $p\leqslant q$ implies $b(p)\subset b(q)$. In other words, $\cF$ is a strong copy of some \textit{extension} of $(P,\leqslant)$. $\LAw(n,P)$ is the largest size of a family $\cF\subseteq 2^{[n]}$ with no weak copies of $P$. As opposed to $\RLA(n,P)$, the problem of determining $\RLAw(n,P)$ is trivial (modulo a result of Erd\H os) for all posets. The family consisting of the middle $|P|-1$ layers of $2^{[n]}$ can be colored by layers without admitting rainbow copies of $P$ as we use only $|P|-1$ colors. On the other hand, the aforementioned result of Erd\H os \cite{E} states that any family $\cF\subseteq 2^{[n]}$ of size $\Sigma(n,|P|-1)+1$ contains a chain of size $|P|$, a weak copy of $P$ which is always rainbow in any proper coloring.

Let us also mention that the literature on the forbidden subposet problem mostly focuses on weak copies (exceptions are \cite{BJ,LM,Me,MP,P,T}), this paper hopefully directs the attention to the less studied version of strong copies.

\medskip

The examples of Secion 2 for posets in $M(A_k)$ have bounded height and width with respect to $k$. Are there any posets $P$ that behave differently (even $A_k$ themselves) or do there exist functions $f$ and $g$ such that if $Q\in M(P)$, then $h(Q)\le f(|P|)$ and $w(Q)\le g(|P|)$ hold? If the answer is positive in both cases, then by Mirsky's or Dilworth's theorem, $|Q|$ is also bounded and it would yield $M(P)$ is finite for all posets $P$.


\begin{thebibliography}{99}
    \bibitem{BJ}
    E. Boehnlein, T. Jiang, Set families with a forbidden induced subposet. Combinatorics, Probability and Computing, 21(4) (2012), 496-511.
    \bibitem{B}
    B. Bukh, Set families with a forbidden subposet. the Electronic Journal of Combinatorics, 16 (2009), R142.
    \bibitem{EIL}
    D. Ellis, M.R. Ivan, I. Leader, Tur\'an densities for daisies and hypercubes. Bulletin of the London Mathematical Society, 56(12) (2024), 3838-3853.
    \bibitem{E}
     P. Erd\H os, On a lemma of Littlewood and Offord, Bull. Amer. Math. Soc. 51 (1945), 898--902
     \bibitem{GP}
     D. Gerbner, B.Patk\'os, \textit{Extremal finite set theory.} (2018) Chapman and Hall/CRC.
     \bibitem{GL}
     J.R. Griggs, W.T. Li, The partition method for poset-free families. J Comb  Optim 25, 587–596 (2013).
     \bibitem{GLsurv}
     J.R. Griggs, W.T. Li, Progress on poset-free families of subsets. In Recent trends in combinatorics (2016) pp. 317-338. Cham: Springer International Publishing.
     \bibitem{GLL}
     J.R. Griggs, W.T. Li, L. Lu, Diamond-free families. Journal of Combinatorial Theory, Series A, 119(2) (2012), 310-322.
     \bibitem{GLu}
     J.R. Griggs, L. Lu, On families of subsets with a forbidden subposet. Combinatorics, Probability and Computing, 18(5) (2009), 731-748. 
     \bibitem{KT}
      G.O.H. Katona, T. Tarj\'an, Extremal problems with excluded subgraphs in the n-cube, Graph Theory, Lagow, 1981, Lecture Notes in Math. 1018 (Springer-Verlag, Berlin, 1983) 084-093.
      \bibitem{LM}
      L. Lu, K.G. Milans, Set families with forbidden subposets, Journal of Combinatorial Theory, Series A 136 (2015) 126-142.
     \bibitem{Lub}
     D. Lubell, A short proof of Sperner’s lemma, J. Combinatorial Theory 1 (1966),
299. 
     \bibitem{Mesh}
     L. D. Meshalkin, A generalization of Sperner’s theorem on the number of subsets
of a finite set, Teor. Verojatnost. i Primenen. 8 (1963), 219–220.
    \bibitem{Me}
    A. M\'eroueh, A LYM inequality for induced posets. Journal of Combinatorial Theory, Series A, 155 (2018) 398-417.
    \bibitem{MP}
    A. Methuku, D. P\'alv\"olgyi, Forbidden hypermatrices imply general bounds on induced forbidden subposet problems. Combinatorics, Probability and Computing, 26(4) (2017), 593-602.
    \bibitem{M}
    L. Mirsky, A dual of Dilworth's decomposition theorem, American Mathematical Monthly, 78 (8) (1971), 876–877.
    \bibitem{P}
    B. Patk\'os, Induced and Non-induced Forbidden Subposet Problems, The Electronic Journal of Combinatorics (2015), P1.30.
    \bibitem{S}
    E. Sperner, Satz \"uber Untermengen einer endlichen Menge, Math Z, 27 (1928), 544--548.
    \bibitem{T}
    I. Tomon, Forbidden induced subposets of given height, Journal of Combinatorial Theory, Series A 161 (2019) 537-562.
    \bibitem{Yam}
    K. Yamamoto, Logarithmic order of free distributive lattice, J. Math. Soc.
Japan 6 (1954), 343–353.
\end{thebibliography}
\end{document}